# A Conformal Mapping Based Fractional Order Approach for Sub-optimal Tuning of PID Controllers with Guaranteed Dominant Pole Placement

Suman Saha[1,2], Saptarshi Das[2,3], Shantanu Das[4], Amitava Gupta[2,3]

1. Drives and Control System Technology Group, CSIR-Central Mechanical Engineering Research Institute (CMERI), Mahatma Gandhi Avenue, Durgapur-713209, India. Email: s_saha@cmeri.res.in

2. Department of Power Engineering, Jadavpur University, Salt Lake Campus, LB-8, Sector 3, Kolkata-700098, India. Email: saptarshi@pe.jusl.ac.in, amitg@pe.jusl.ac.in

3. School of Nuclear Studies & Applications, Jadavpur University, Salt Lake Campus, LB-8, Sector 3, Kolkata-700098, India.

4. Reactor Control Division, Bhabha Atomic Research Centre, Mumbai-400085, India.
Email: shantanu@barc.gov.in

**ABSTRACT:**
A novel conformal mapping based Fractional Order (FO) methodology is developed in this paper for tuning existing classical (Integer Order) Proportional Integral Derivative (PID) controllers especially for sluggish and oscillatory second order systems. The conventional pole placement tuning via Linear Quadratic Regulator (LQR) method is extended for open loop oscillatory systems as well. The locations of the open loop zeros of a fractional order PID (FOPID or PI$^\lambda$D$^\mu$) controller have been approximated in this paper vis-à-vis a LQR tuned conventional integer order PID controller, to achieve equivalent integer order PID control system. This approach eases the implementation of analog/digital realization of a FOPID controller with its integer order counterpart along with the advantages of fractional order controller preserved. It is shown here in the paper that decrease in the integro-differential operators of the FOPID/PI$^\lambda$D$^\mu$ controller pushes the open loop zeros of the equivalent PID controller towards greater damping regions which gives a trajectory of the controller zeros and dominant closed loop poles. This trajectory is termed as "M-curve". This phenomena is used to design a two-stage tuning algorithm which reduces the existing PID controller's effort in a significant manner compared to that with a single stage LQR based pole placement method at a desired closed loop damping and frequency.

**Nomenclature:**
$\zeta^{ol}$ : Damping ratio of the second order (open loop) plant

$\omega_n^{ol}$ : Natural frequency of the second order (open loop) plant

$K$ : Gain of the second order (open loop) plant (dc gain is $K\big/\left(\omega_n^{ol}\right)^2$ )

$K_p$ : Proportional gain of the PID controller

$K_i$ : Integral gain of the PID controller

$K_d$ : Derivative gain of the PID controller

$\zeta^{cl}$ : Damping ratio of the (second order approximated) closed loop system

$\omega_n^{cl}$ : Damping ratio of the (second order approximated) closed loop system

$Q$ : Symmetric positive semi-definite weighting matrix for LQR cost function

$R$ : Positive weighting factor for LQR cost function

$P$ : Symmetric positive definite solution matrix of the algebraic Riccati equation

$F$ : State-feedback gain matrix of the optimal regulator

$r(t)$ : Set-point

$e(t)$ : Error signal

$u(t)$ : Control signal

$y(t)$ : Response of the controlled system

## 1. Introduction

In most process control applications, dominant pole placement tuning is a popular technique with second order approximations for sluggish or oscillatory processes [1], [2]. It is shown in Wang *et al.* [3] that guaranteed dominant pole placement can be done with PID controllers if real part of the resulting non-dominant closed loop poles is at least 3-5 times larger than that of the dominant closed loop poles. He *et al.* [4] tried to tune PID controllers using LQR based technique taking pole placement into consideration while selecting the weighting matrices for the optimal quadratic regulator. The technique, proposed in [4] ensures pole placement with PID controllers while also minimizing the state deviations and controller effort for doing so [5]. The concept of LQR based PI/PID tuning was primarily derived for first order systems with PI controllers in [4]. Then it was extended for second order systems with PID controllers as well; in order to cancel the open loop system pole with one real controller zero [4]. This reduces the formulation to a simple PI controller design for handling first order systems. This approach fails to give satisfactory results especially for systems having highly oscillatory open loop dynamics i.e. complex conjugate open loop poles. The present methodology developed in this paper, primarily tries to focus on the unresolved issues addressed in [3], [4] e.g.

(a) Condition for guaranteed dominant pole placement

(b) Extension of LQR based PID tuning for highly oscillatory and sluggish processes with the same pole placement technique

(d) Comparison of the control cost and also initial controller efforts from the regulator design point of view.

FOPID or $PI^\lambda D^\mu$ controllers [6] recently have become popular in process control as it has two extra parameters for tuning i.e. the integro-differential orders; giving us extra freedom of tuning. Therefore, FOPID controller is expected to give better performance over conventional PID controllers [6], [7]. While conformal mapping for FOPID is done with *s* to *w* plane-transform, the transfer function of a FOPID controller [6], [7] places two fractional zeros and one fractional pole at the origin in the complex *w*-plane; with differ-integral fractional orders of the controller are equal. The solution of fractional order systems can run into multiple Riemann Sheets [7]. Let us assume that the domain of existence of the controller zeros and poles are limited to only the primary Riemann-sheet, so as to ensure that the closed loop system is not hyper-damped or ultra-damped which does not come into dominant dynamics [7]-[11]; then position of the two fractional order zeros in complex *s*-plane can always be replaced by two integer order zeros. This is because guaranteed pole placement is feasible with a PID controller anywhere in the negative real part of the *s*-plane. With this approach, our focus is to represent the approximated PID controller gains in terms of FOPID controller parameters which reduce the complicacy of FO controller realization; while also preserving the time response of the controlled system.

In this paper, it is observed that with the decrease in the order of a FOPID controller, the position of the zeros shifts towards greater damping regions and then move very fast towards smaller damping and ultimately goes towards instability. While the obtained trajectory of the controller zeros or the dominant closed loop poles takes a certain pattern and we termed as "M-curve". The proposed methodology first assumes a LQR based PID controller with desired location of the dominant pole placement. Also, it is always feasible to place the dominant closed loop poles at the same location of the primary Riemann sheet with a FOPID controller. Now with variation in the FO controller orders, the equivalent integer order PID controller with fractional zeros in the same location gives different set of controller gains. Upon approximation with integer order PID, the controller does not remain optimal but is suboptimal. This conformal mapping based sub-optimal controller design possesses certain strengths over the LQR based optimal PID controller, which is illustrated and elucidated with credible and exhaustive numerical simulations.

The objective of this paper is to put forward a novel methodology that tunes a PID controller with an LQR based dominant pole placement method at a lower damping than the desired one in the first stage; and then considering FOPID controller zeros at the same location for pole-placement. Thereafter the order of FO controller is decreased so as to obtain the approximated integer order suboptimal PID gains that forces the closed loop poles to move towards greater (desired) damping. Simulation results are given to justify that the proposed two-stage tuning of a PID controller significantly reduces the control signal vis-à-vis a single stage LQR based PID controller to achieve the same desired closed loop damping that is percentage overshoot.

The rest of the paper is organized as follows. Section 2 proposes LQR based guaranteed dominant pole placement of integer order PID controller with highly oscillatory processes too [12] and discusses about the inverse optimal control costs involved in pole placement problem. Section 3 introduces a new fractional order approach of PID controller tuning. Simulation studies are carried out in section 4 for three different class of second order processes to show the effect of variation in orders of FO controller and its effect on the integer order approximation of the fractional order controller (zeros and dominant closed loop poles). The advantage of two-stage sub-optimal tuning of PID controllers over the conventional LQR based optimal tuning is dealt in section 5. The paper ends with the conclusion as section 6, followed by the references.

## 2. Concept of LQR based guaranteed dominant pole placement with PID controllers for second order systems

### 2.1. Criteria for guaranteed dominant pole placement

In this section, a brief idea is presented regarding the accuracy of guaranteed pole placement with PID controllers. Let us consider that a second order process $G_p$ (with sluggish "$S$" shaped or oscillatory open loop dynamics) needs to be controlled with an integer order PID controller $C(s)$ of the form (1).

$$C(s) = K_p + \frac{K_i}{s} + K_d s = \frac{K_d s^2 + K_p s + K_i}{s} \tag{1}$$

Here, the process is characterized by the open loop transfer function

$$G_p(s) = \frac{K}{s^2 + 2\zeta^{ol}\omega_n^{ol}s + \left(\omega_n^{ol}\right)^2} \tag{2}$$

Then, the closed loop transfer function becomes

$$
\begin{aligned}
G_{cl} &= \frac{G_p C}{1 + G_p C} \\
&= \frac{K\left(K_d s^2 + K_p s + K_i\right)}{s^3 + s^2\left(2\xi^{ol}\omega_n^{ol} + KK_d\right) + s\left(\left(\omega_n^{ol}\right)^2 + KK_p\right) + KK_i}
\end{aligned}
\tag{3}
$$

From (2) it is clear that the open loop plant has two poles at $\left[-\zeta^{ol}\omega_n^{ol} \pm j\omega_n^{ol}\sqrt{1 - \left(\zeta^{ol}\right)^2}\right]$ and from (1) is also evident that the PID controller has one pole at origin and also two zeros at $\left[-\frac{K_p}{2K_d} \pm j\sqrt{\frac{K_i}{K_d} - \frac{K_p^2}{4K_d^2}}\right]$. It is considered that both the process poles and controller zeros are complex conjugates in the complex $s$-plane. From (3), it is seen that the closed loop system has two zeros and three poles in the complex $s$-plane and position of the closed loop zeros remains unchanged as in (1) while position of the closed loop poles changes depending on the PID controller gains. For guaranteed pole placement with PID controllers, the closed loop system (3) should have one real pole which should be far away from the real part of the other two complex (conjugate) closed loop poles as discussed in Wang *et al.* [3]. It is also reported in [3] that the contribution of the real pole in the closed loop dynamics becomes insignificant if magnitude of the real closed loop pole be at least 3-5 (let us call this factor as $m$ 'relative dominance') times greater than real part of the complex closed loop poles. Now, if the desired closed loop performance of a second order system be given as the specifications on $\zeta^{cl}$ and $\omega_n^{cl}$ as in [4], one can easily replace the position of the real zero ($\alpha$) by $\left(-m\zeta^{cl}\omega_n^{cl}\right)$, provided $\alpha = -m\zeta^{cl}\omega_n^{cl}$ is chosen to be large enough with respect to $\left(-\zeta^{cl}\omega_n^{cl}\right)$; while gradually increasing the value of the relative dominance ($m$) and checking the accuracy of pole placement in time as well as frequency domain. Otherwise, significant effect of the real pole will be there on the closed loop dynamics manifesting as sluggishness and therefore dominant pole placement can not be guaranteed. With a proper choice of $m$, the third order closed loop system (3) will behave like an almost second order system having the user specified damping $\zeta^{cl}$ (percentage of maximum overshoot) and frequency $\omega_n^{cl}$ (rise time). In this circumstances, the characteristic polynomial for $G_{cl}$ is written as:

$$\left(s + m\zeta^{cl}\omega_n^{cl}\right)\left(s^2 + 2\zeta^{cl}\omega_n^{cl}s + \left(\omega_n^{cl}\right)^2\right) = 0$$
$$\Rightarrow \left[s^3 + s^2\left(2+m\right)\zeta^{cl}\omega_n^{cl} + s\left(\left(\omega_n^{cl}\right)^2 + 2m\left(\zeta^{cl}\right)^2\left(\omega_n^{cl}\right)^2\right) + m\zeta^{cl}\left(\omega_n^{cl}\right)^3\right] = 0 \tag{4}$$

Comparing the coefficients of (4) with the denominator of (3), the PID controller parameters are obtained .Therefore,

$$\left.\begin{array}{l} K_p = \dfrac{\left(\omega_n^{cl}\right)^2 + 2m\left(\zeta^{cl}\right)^2\left(\omega_n^{cl}\right)^2 - \left(\omega_n^{ol}\right)^2}{K} \\[4mm] K_i = \dfrac{m\zeta^{cl}\left(\omega_n^{cl}\right)^3}{K} \\[4mm] K_d = \dfrac{\left(2+m\right)\zeta^{cl}\omega_n^{cl} - 2\zeta^{ol}\omega_n^{ol}}{K} \end{array}\right\} \tag{5}$$

To choose a suitable value of $m$ let us take a highly oscillatory system excluding the delay as reported in Panda, Yu and Huang [13]. Here, the open loop system parameters are $K = 1, \zeta^{ol} = 0.2, \omega_n^{ol} = 0.1\,\text{rad/sec}$ and the desired parameters of closed loop system are $\zeta^{cl} = 0.98, \omega_n^{cl} = 2$ rad/sec. Using the relations in (5), the gains of the PID controller (1) are calculated which produces exact pole placement at desired damping and frequency, provided $m$ is chosen judiciously (by iteratively checking the accuracy of pole placement). The corresponding root locus and unit step responses are shown in Fig.1.

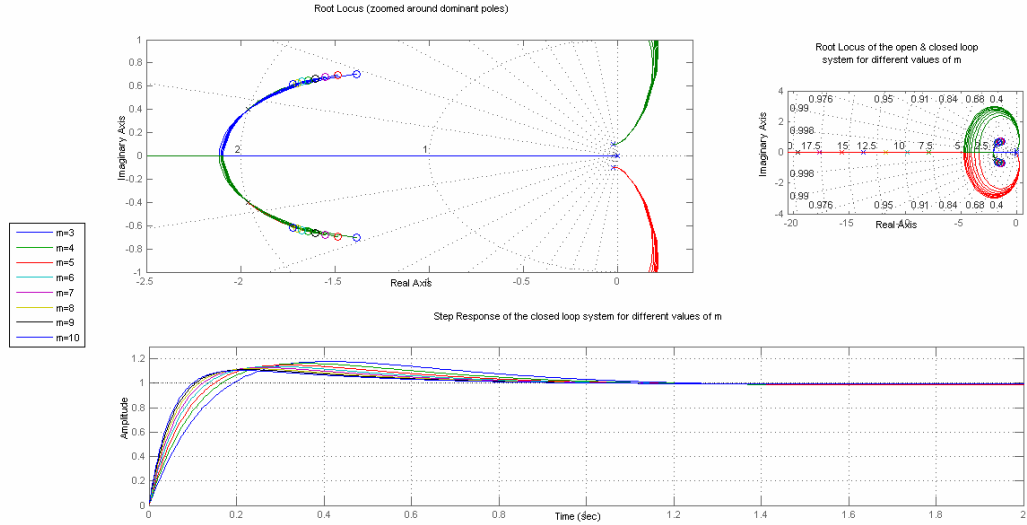

Fig. 1: Root locus and unit step response for increasing $m$.

From the root locus plots in Fig. 1, it is evident that the method places the closed loop poles exactly at a desired damping and frequency. But a small reduction in maximum overshoot ($\%M_p$) and rise time ($t_r$) can be noticed with gradual increase in $m$. It is also observed that further increase of 'relative dominance' above $m = 10$ does not affect the time response in a considerable manner even for a highly oscillatory process. For this reason of saturation in $m$ in this paper, we choose $m = 10$ for all cases to find out the controller parameters with (5) in the subsequent sections.

### 2.2. Design of LQR based optimal PID controller for second order systems:

He *et al.* [4] have given a formulation for tuning over-damped and critically-damped second order systems which has been extended here for lightly damped processes as well. Additionally, in [4], it is suggested that one of the real poles can be cancelled out by placing one of the controller zeros at the same position on the negative real axis in complex $s$-plane. Thus the second order plant to be controlled with a PID controller can be reduced to a first order process to be controlled by a PI controller. Indeed, this approach of He *et al.* [4] is not valid for lightly damped processes since a single complex pole of the process cannot be eliminated with a single complex zero of controller, as they are in conjugate pairs. With this approach of optimal PID tuning for second order processes of [4], the provision of simultaneously and optimally finding the three parameters of a PID controller (i.e. $K_p, K_i, K_d$) is lost which is addressed in this paper. The present approach takes the error, error rate and integral of error as the state variables and designs the optimal controller gains as the parameters of the PID regulator (Fig. 2).

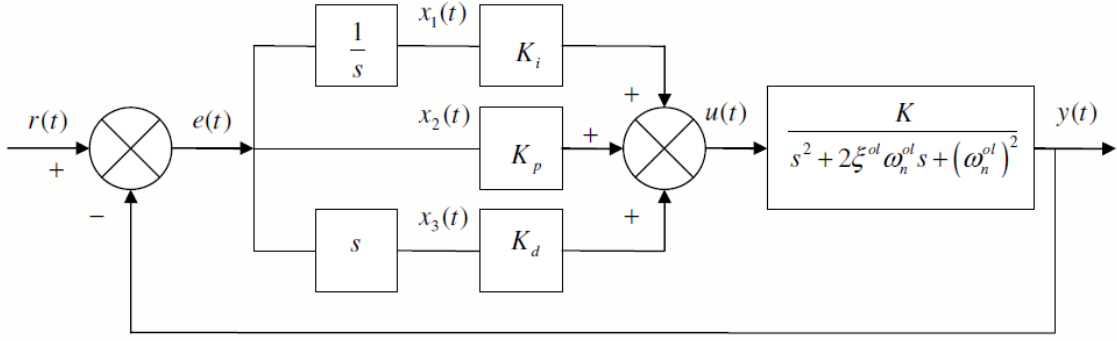

Fig. 2. Formulation of LQR based PID controller for second order processes.

Therefore, the state variables are:

$$x_1 = \int e(t)dt, \quad x_2 = e(t), \quad x_3 = \frac{de(t)}{dt} \tag{6}$$

From the block diagram presented in Fig. 2, we obtain

$$\frac{Y(s)}{U(s)} = \frac{K}{s^2 + 2\xi^{ol}\omega_n^{ol}s + \left(\omega_n^{ol}\right)^2} = \frac{-E(s)}{U(s)} \tag{7}$$

In the case of state feedback regulator, the external set-point does not affect the controller design therefore $r(t) = 0$, in Fig. 2. In (7), the relation $y(t) = -e(t)$ is valid for standard regulator problem as also reported in He *et al.* [4], when there is no change in the set point. Thus, equation (7) takes the form

$$\left[s^2 + 2\xi^{ol}\omega_n^{ol}s + \left(\omega_n^{ol}\right)^2\right]E(s) = -KU(s) \tag{8}$$

$$\Rightarrow \quad \ddot{e} + 2\xi^{ol}\omega_n^{ol}\dot{e} + \left(\omega_n^{ol}\right)^2 e = -Ku \tag{9}$$

Using (6), equation (9) can be re-written as

$$\dot{x}_3 + 2\xi^{ol}\omega_n^{ol}x_3 + \left(\omega_n^{ol}\right)^2 x_2 = -Ku \tag{10}$$

Using (6) and (10) the state space formulation of the above system becomes

$$\begin{bmatrix} \dot{x}_1 \\ \dot{x}_2 \\ \dot{x}_3 \end{bmatrix} = \begin{bmatrix} 0 & 1 & 0 \\ 0 & 0 & 1 \\ 0 & -\left(\omega_n^{ol}\right)^2 & -2\xi^{ol}\omega_n^{ol} \end{bmatrix} \begin{bmatrix} x_1 \\ x_2 \\ x_3 \end{bmatrix} + \begin{bmatrix} 0 \\ 0 \\ -K \end{bmatrix} u \tag{11}$$

Comparing (11) with the standard state-space representation of a system, that is (12)

$$\dot{x}(t) = Ax(t) + Bu(t) \tag{12}$$

we get the system matrices as (13)

$$A = \begin{bmatrix} 0 & 1 & 0 \\ 0 & 0 & 1 \\ 0 & -\left(\omega_n^{ol}\right)^2 & -2\xi^{ol}\omega_n^{ol} \end{bmatrix}, \quad B = \begin{bmatrix} 0 \\ 0 \\ -K \end{bmatrix} \tag{13}$$

In order to obtain an LQR formulation with the system (13), the following quadratic cost should be minimized, that is (14)

$$J = \int_0^\infty \left[ x^T(t)Qx(t) + u^T(t)Ru(t) \right] dt \qquad (14)$$

It is shown in [5] that minimization of (14) gives the state feedback control signal as

$$u(t) = -R^{-1}B^T Px(t) = -Fx(t) \qquad (15)$$

where, $P$ is the symmetric positive definite solution of the Continuous Algebraic Riccati Equation (CARE) given by (16)

$$A^T P + PA - PBR^{-1}B^T P + Q = 0 \qquad (16)$$

Here, the weighting matrix $Q = \begin{bmatrix} Q_1 & 0 & 0 \\ 0 & Q_2 & 0 \\ 0 & 0 & Q_3 \end{bmatrix}$ is symmetric positive semi-definite and the weighting factor $R$ is a positive number. It is a common practice in optimal control to design regulators with varying $Q$, while keeping $R$ fixed [14] and generally they are designed with user specified closed loop performance specifications [4]. Considering that the unique solution of the CARE (16) as $P = \begin{bmatrix} P_{11} & P_{12} & P_{13} \\ P_{12} & P_{22} & P_{23} \\ P_{13} & P_{23} & P_{33} \end{bmatrix}$, the state feedback gain

matrix corresponding to the optimal control signal in (15) therefore is

$$F = R^{-1}B^T P = R^{-1} \begin{bmatrix} 0 & 0 & -K \end{bmatrix} \begin{bmatrix} P_{11} & P_{12} & P_{13} \\ P_{12} & P_{22} & P_{23} \\ P_{13} & P_{23} & P_{33} \end{bmatrix}$$

$$= -R^{-1}K \begin{bmatrix} P_{13} & P_{23} & P_{33} \end{bmatrix} = -\begin{bmatrix} K_i & K_p & K_d \end{bmatrix} \qquad (17)$$

Using (6), the corresponding expression for the control signal is obtained as (18)

$$u(t) = -Fx(t) = -\begin{bmatrix} -K_i & -K_p & -K_d \end{bmatrix} \begin{bmatrix} x_1(t) \\ x_2(t) \\ x_3(t) \end{bmatrix}$$

$$= K_i \int e(t)dt + K_p e(t) + K_d \frac{de(t)}{dt} \qquad (18)$$

From (17), the third row (or column) of symmetric positive definite matrix $P$ is obtained in terms of PID controller gains as in (19)

$$P_{13} = \frac{K_i}{R^{-1}K}, \quad P_{23} = \frac{K_p}{R^{-1}K}, \quad P_{33} = \frac{K_d}{R^{-1}K} \qquad (19)$$

Now the closed loop system matrix for the system (13) with state feedback gains (17) is as follows

$A_c = A - BF$

$$= \begin{bmatrix} 0 & 1 & 0 \\ 0 & 0 & 1 \\ 0 & -\left(\omega_n^{ol}\right)^2 & -2\xi^{ol}\omega_n^{ol} \end{bmatrix} + \begin{bmatrix} 0 \\ 0 \\ -K \end{bmatrix} R^{-1}K \begin{bmatrix} P_{13} & P_{23} & P_{33} \end{bmatrix} \quad (20)$$

$$= \begin{bmatrix} 0 & 1 & 0 \\ 0 & 0 & 1 \\ \left(-R^{-1}K^2 P_{13}\right) & \left(-\left(\omega_n^{ol}\right)^2 - R^{-1}K^2 P_{23}\right) & \left(-2\xi^{ol}\omega_n^{ol} - R^{-1}K^2 P_{33}\right) \end{bmatrix}$$

The corresponding characteristics equation for the closed loop system is (21)

$\Delta(s) = |sI - A_c| = 0$

$$\Rightarrow \quad s^3 + s^2 \left(2\xi^{ol}\omega_n^{ol} + R^{-1}K^2 P_{33}\right) + s\left(\left(\omega_n^{ol}\right)^2 + R^{-1}K^2 P_{23}\right) + R^{-1}K^2 P_{13} = 0 \quad (21)$$

Now comparing the coefficients of (21) with the denominator of (3) we get,

$$\left.\begin{aligned} 2\xi^{ol}\omega_n^{ol} + R^{-1}K^2 P_{33} &= \left(2+m\right)\zeta^{cl}\omega_n^{cl} \\ \left(\omega_n^{ol}\right)^2 + R^{-1}K^2 P_{23} &= \left(\omega_n^{cl}\right)^2 + 2m\left(\zeta^{cl}\right)^2 \left(\omega_n^{cl}\right)^2 \\ R^{-1}K^2 P_{13} &= m\zeta^{cl}\left(\omega_n^{cl}\right)^3 \end{aligned}\right\} \quad (22)$$

From (22), by knowing the open loop process characteristics (i.e. $\xi^{ol}, \omega_n^{ol}, K$) and the desired closed loop system dynamics ($\zeta^{cl}, \omega_n^{cl}, m$) the elements of the third row of matrix $P$ is solved as in (23)

$$\left.\begin{aligned} P_{13} &= \frac{m\zeta^{cl}\left(\omega_n^{cl}\right)^3}{R^{-1}K^2} \\ P_{23} &= \frac{\left(\omega_n^{cl}\right)^2 + 2m\left(\zeta^{cl}\right)^2\left(\omega_n^{cl}\right)^2 - \left(\omega_n^{ol}\right)^2}{R^{-1}K^2} \\ P_{33} &= \frac{\left(2+m\right)\zeta^{cl}\omega_n^{cl} - 2\xi^{ol}\omega_n^{ol}}{R^{-1}K^2} \end{aligned}\right\} \quad (23)$$

Now, in order to solve the Continuous Algebraic Riccati Equation CARE (16), the following set of linear equations (24) is to be solved

$$\begin{bmatrix} 0 & 1 & 0 \\ 0 & 0 & 1 \\ 0 & -\left(\omega_n^{ol}\right)^2 & -2\xi^{ol}\omega_n^{ol} \end{bmatrix}\begin{bmatrix} P_{11} & P_{12} & P_{13} \\ P_{12} & P_{22} & P_{23} \\ P_{13} & P_{23} & P_{33} \end{bmatrix} + \begin{bmatrix} P_{11} & P_{12} & P_{13} \\ P_{12} & P_{22} & P_{23} \\ P_{13} & P_{23} & P_{33} \end{bmatrix}\begin{bmatrix} 0 & 1 & 0 \\ 0 & 0 & 1 \\ 0 & -\left(\omega_n^{ol}\right)^2 & -2\xi^{ol}\omega_n^{ol} \end{bmatrix}$$

$$-\begin{bmatrix} P_{11} & P_{12} & P_{13} \\ P_{12} & P_{22} & P_{23} \\ P_{13} & P_{23} & P_{33} \end{bmatrix}\begin{bmatrix} 0 \\ 0 \\ -K \end{bmatrix}R^{-1}\begin{bmatrix} 0 & 0 & -K \end{bmatrix}\begin{bmatrix} P_{11} & P_{12} & P_{13} \\ P_{12} & P_{22} & P_{23} \\ P_{13} & P_{23} & P_{33} \end{bmatrix} + \begin{bmatrix} Q_1 & 0 & 0 \\ 0 & Q_2 & 0 \\ 0 & 0 & Q_3 \end{bmatrix} = 0 \qquad (24)$$

Now from (23), the third row or third column of symmetric positive definite $P$ matrix (solution of Riccati equation) is known from the open loop and desired closed loop system dynamics. With these known parameters using (24), the other elements of $P$ and $Q$ matrix can then be evaluated as following (25) and (26)

$$\left.\begin{aligned} P_{11} &= \left(\omega_n^{ol}\right)^2 P_{13} + R^{-1}K^2 P_{13}P_{23} \\ P_{12} &= 2\xi^{ol}\omega_n^{ol}P_{13} + R^{-1}K^2 P_{13}P_{33} \\ P_{22} &= 2\xi^{ol}\omega_n^{ol}P_{23} + R^{-1}K^2 P_{23}P_{33} + \left(\omega_n^{ol}\right)^2 P_{33} - P_{13} \end{aligned}\right\} \qquad (25)$$

$$\left.\begin{aligned} Q_1 &= R^{-1}K^2 P_{13}^2 \\ Q_2 &= R^{-1}K^2 P_{23}^2 - 2\left(P_{12} - \left(\omega_n^{ol}\right)^2 P_{23}\right) \\ Q_3 &= R^{-1}K^2 P_{33}^2 - 2\left(P_{23} - 2\xi^{ol}\omega_n^{ol}P_{33}\right) \end{aligned}\right\} \qquad (26)$$

Next an example is shown to illustrate the equivalent inverse optimal control formulation for such a pole placement problem. For the system defined in section 2.1, with parameters $K = 1, \zeta^{ol} = 0.2, \omega_n^{ol} = 0.1$ rad/sec and desired parameters of closed loop system defined as $\zeta^{cl} = 0.98, \omega_n^{cl} = 2$ rad/sec, $m = 10$, the weighting matrices ($Q$ and $R$) and the solution matrix of CARE ($P$) are respectively calculated from (26), (23) and (25) and reported in (27) (28) and (29)

$$Q_H = \begin{bmatrix} 6.1466 & 0 & 0 \\ 0 & 2.8459 & 0 \\ 0 & 0 & 0.3915 \end{bmatrix} \times 10^3, \quad R_H = 1 \qquad (27)$$

$$P_H = \begin{bmatrix} 6.3372 & 1.8440 & 0.0784 \\ 1.8440 & 1.8228 & 0.0808 \\ 0.0784 & 0.0808 & 0.0235 \end{bmatrix} \times 10^3 \qquad (28)$$

The system matrices are $A = \begin{bmatrix} 0 & 1 & 0 \\ 0 & 0 & 1 \\ 0 & -0.01 & -0.04 \end{bmatrix}, B = \begin{bmatrix} 0 \\ 0 \\ -1 \end{bmatrix} \qquad (29)$

The corresponding state feedback gain matrix is obtained as (30)

$$F_H = -\begin{bmatrix} 78.400 & 80.822 & 23.480 \end{bmatrix}$$
$$= -\begin{bmatrix} K_i^H & K_p^H & K_d^H \end{bmatrix} \tag{30}$$

Now, the dominant pole placement could have been formulated for placing the closed loop poles at a relatively smaller damping. Here we investigate and note the difference in control cost. In order to do so, let us consider the closed loop specification of the same plant with parameters as $\zeta^{cl} = 0.75$, $\omega_n^{cl} = 2$ rad/sec, $m = 10$. The corresponding weighting matrices and solution matrix for CARE are (31), (32)

$$Q_L = \begin{bmatrix} 3.600 & 0 & 0 \\ 0 & 0.2410 & 0 \\ 0 & 0 & 0.2260 \end{bmatrix} \times 10^3, \quad R_L = 1 \tag{31}$$

$$P_L = \begin{bmatrix} 2.940 & 1.080 & 0.060 \\ 1.080 & 0.822 & 0.049 \\ 0.060 & 0.049 & 0.018 \end{bmatrix} \times 10^3 \tag{32}$$

The corresponding state feedback gain matrix is obtained as

$$F = -\begin{bmatrix} 60.00 & 48.99 & 17.96 \end{bmatrix}$$
$$= -\begin{bmatrix} K_i^L & K_p^L & K_d^L \end{bmatrix} \tag{33}$$

In equations (27)-(33), subscripts of matrices and superscripts of controller gains ($H$ and $L$) corresponds to the design at higher and lower closed loop damping respectively. Also, it is shown in [15] that the value of the integral performance index (14) can be calculated from the Riccati solution ($P$ matrix) using the initial values of the state variables i.e.

$$J = x^T(0)Px(0) = \int_0^\infty \left[ x^T(t)Qx(t) + u^T(t)Ru(t) \right] dt \tag{34}$$

For a PID controller as in our case, initial values of the state variables (i.e. error, error rate and integral of error) cannot be calculated directly to find out the optimal control cost (34). The reason is since with a step-input excitation the initial value of the error rate will tend to infinity and initial value of integral error will tend to zero with the initial value of error signal remaining unity. To overcome this problem the following methodology is adopted. From the respective Riccati solutions (28) and (32), we get the eigen-values of the differential Riccati solutions as (35)

$$eig\begin{bmatrix} P_H - P_L \end{bmatrix} = eig\begin{bmatrix} \begin{pmatrix} 3.3972 & 0.7640 & 0.0184 \\ 0.7640 & 1.0008 & 0.0318 \\ 0.0184 & 0.0318 & 0.0055 \end{pmatrix} \times 10^3 \end{bmatrix} = \begin{bmatrix} 0.0045 \\ 0.7788 \\ 3.6203 \end{bmatrix} \times 10^3 \tag{35}$$

Clearly the eigen-values of $\begin{bmatrix} P_H - P_L \end{bmatrix}$ are positive which indicates that matrix $\begin{bmatrix} P_H - P_L \end{bmatrix}$ is positive definite. Now for all initial value of the state variables i.e. $x(0)$, pre-multiplication with $x^T(0)$ and post-multiplication with $x(0)$, the Riccati solution matrices for high and low closed loop damping take the form (36), giving comparison of the achievable cost of control.

$$x^T(0)P_H x(0) > x^T(0)P_L x(0) \qquad \Rightarrow J_H > J_L \tag{36}$$

Equation (36) implies that cost of pole placement for greater damping with PID controllers via LQR is more than that for smaller closed loop damping.

## 3. Fractional zero placement approach for FOPID controllers and conformal mapping based approximation in complex $s \leftrightarrow w$ plane:

Let us consider a FOPID or $PI^{\lambda}D^{\mu}$ controller (37) with integral order ($\lambda$) and the derivative order ($\mu$) set to the same value ($q$), i.e. $\lambda = \mu = q$ similar to that in [16]-[17]. Therefore, the transfer function of the FOPID controller becomes

$$\widetilde{C}(s) = K_p + \frac{K_i}{s^{\lambda}} + K_d s^{\mu} = K_p + \frac{K_i}{s^q} + K_d s^q = \frac{K_d s^{2q} + K_p s^q + K_i}{s^q} \qquad (37)$$

First we shall transform the controller (37) to $w$-plane, by putting variable transformation as $s^q = w$. This is called conformal mapping [7], and the transfer function is now in $w$-plane. The controller is transformed as; $\widetilde{C}(w) = (K_d w^2 + K_p w + K_i)/w$ in $w$-plane. Clearly, the controller $\widetilde{C}(w)$ places two fractional zeros and one fractional pole in the complex $w$ plane. Now, the position of the fractional zeros in complex $w$ plane can be determined by solving the following equation:

$$K_d s^{2q} + K_p s^q + K_i = 0 \qquad \Rightarrow K_d w^2 + K_p w + K_i = 0 \qquad (38)$$

$$\Rightarrow w_{1,2} = \frac{-K_p \pm \sqrt{K_p^2 - 4K_i K_d}}{2K_d} \qquad (39)$$

Now if $K_p^2 > 4K_i K_d$, the two controller zeros in $w$ plane becomes real (ultra-damped), which is not desirable since roots of only primary Riemann sheet contributes to the dominant dynamics of a system [7], [9]. If $K_p^2 < 4K_i K_d$, the two zeros of the FOPID controller (37) becomes complex conjugates as mentioned in (40)

$$w_{1,2} = \frac{-K_p}{2K_d} \pm j \frac{\sqrt{4K_i K_d - K_p^2}}{2K_d} \qquad (40)$$

If $\varphi$ be the angle, made by the line joining the origin and controller zeros with the positive real axis as shown in Fig. 3, then (40) can be rewritten as

$$w_{1,2} = r\left(\cos\varphi + j\sin\varphi\right) = \left(s_{1,2}\right)^q \qquad (41)$$

Comparing (40) and (41) we get,

$$r\cos\varphi = \frac{-K_p}{2K_d}, \quad r\sin\varphi = \frac{\sqrt{4K_i K_d - K_p^2}}{2K_d} \qquad (42)$$

$$\Rightarrow \varphi = \pi - \tan^{-1}\left(\frac{\sqrt{4K_i K_d - K_p^2}}{K_p}\right), \quad r = \sqrt{\frac{K_i}{K_d}} \qquad (43)$$

Now, position of the FOPID controller zeros as complex conjugates in complex $w$ plane do not guarantee that they will not contribute to unstable ($\varphi < \frac{\pi q}{2}$), hyper-damped ($\varphi > \pi q$) or ultra-damped dynamics ($\varphi = \pi$) [7]-[8]. In order to ensure that the controller zeros are placed in the under-damped region in complex $w$ plane (and corresponding negative real part of the complex $s$ plane), the following criteria must be satisfied:

$$\frac{\pi q}{2} < \varphi < \pi q \tag{44}$$

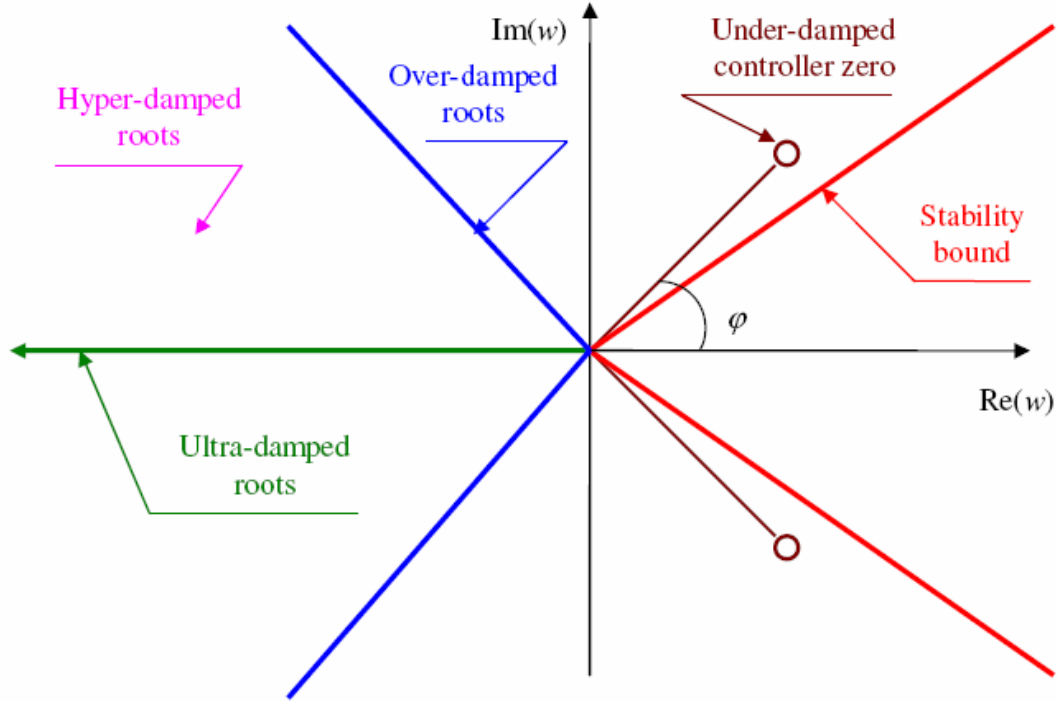

Fig. 3. Position of the controller zeros in complex $w$ plane.

Therefore using conformal mapping based approach the under-damped (complex conjugate) zeros of the FOPID controller in $w$ plane can be mapped back in the complex $s$ plane as

$$s_{1,2} = r^{1/q} \left( \cos \varphi \pm j \sin \varphi \right)^{1/q} = r^{1/q} \left( \cos \frac{\varphi}{q} \pm j \sin \frac{\varphi}{q} \right) \tag{45}$$

Using (43) from (45) we get

$$s_{1,2} = \left( \left( \frac{K_i}{K_d} \right)^{1/2q} \cos \left( \frac{\varphi}{q} \right) \right) \pm j \left( \left( \frac{K_i}{K_d} \right)^{1/2q} \sin \left( \frac{\varphi}{q} \right) \right) \tag{46}$$

Clearly, (46) represents two conjugate fractional ($q^{th}$) order zeros which can be approximated using two integer ($1^{st}$) order zeros, since PID controller is capable of placing zeros anywhere in the negative real part of the $s$ plane with suitable modification of its gains. The two zeros are to be placed at the same locations as (46) with a suitable choice of the fractional order $q$. Therefore, numerator of controller $C(s)$ in (1) with replacement of the fractional order zeros of controller $\widetilde{C}(s)$ in (37) with integer order zeros yields

$$(s - s_1)(s - s_2) = 0 \qquad \Rightarrow s^2 - s(s_1 + s_2) + s_1 s_2 = 0 \tag{47}$$

From (46) the sum and product of the roots can be calculated as

$$\left(s_1 + s_2\right) = 2\left(\frac{K_i}{K_d}\right)^{1/2q}\left(\cos\frac{\theta}{q}\right), \quad s_1 s_2 = \left(\frac{K_i}{K_d}\right)^{1/q} \tag{48}$$

Now from equation (47) and (48) we get

$$s^2 - 2\left(\frac{K_i}{K_d}\right)^{1/2q}\left(\cos\frac{\varphi}{q}\right)s + \left(\frac{K_i}{K_d}\right)^{1/q} = 0$$

$$\Rightarrow \left[\left(K_d\right)^{1/q}\right]s^2 + \left[-2\left(K_i K_d\right)^{1/2q}\left(\cos\frac{\varphi}{q}\right)\right]s + \left[\left(K_i\right)^{1/q}\right] = 0 \tag{49}$$

$$\Rightarrow \widehat{K_d}s^2 + \widehat{K_p}s + \widehat{K_i} = 0$$

With the above formulation (49) a $PI^{\lambda}D^{\mu}$ controller with $\lambda = \mu = q$ can be replaced by a simple PID controller whose zero positions changes with variation in $q$. Since, hardware implementation of fractional order elements and FOPID controllers are difficult due its infinite dimensional nature, therefore the true potential of a FOPID controller [16] can be implemented in a relatively simpler manner with an equivalent PID controller. The parameters of the approximated PID controller are

$$\left.\begin{aligned}\widehat{K_d} &= \left(K_d\right)^{1/q} \\ \widehat{K_p} &= -2\left(K_i K_d\right)^{1/2q}\left(\cos\frac{\varphi}{q}\right) \\ \widehat{K_i} &= \left(K_i\right)^{1/q}\end{aligned}\right\} \tag{50}$$

In (50) the expression for the proportional gain includes a negative sign which is counter intuitive! This can be justified in a way that with the approximated PID controller (with parameters $\widehat{K_p}, \widehat{K_i}, \widehat{K_d}$) for FOPID controller (37), the stability condition of under damped region in $w$-plane (44), has the value of $\cos(\varphi/q)$, which is negative in this region of operation, that is $(\pi q/2) < \varphi < \pi q$ in Fig. 3, giving positive values for controller gain $\widehat{K_p}$ in (50).

## 4. Effect of variation in controller's differ-integral orders on the closed loop control performance

The closed loop performance of a well tuned FOPID controller gets heavily affected by slight change in its differ-integral orders [7]. The variation in the order of the FOPID controller (37) in the pre defined way leads to an advanced tuning strategy for the approximated PID controller with gains (50) which is described in the next section for three different classes of second order systems.

### 4.1. Under-damped process

Panda, Yu and Huang [13] have shown the tuning results of various second order processes. In this paper, the following (51) lightly damped process is considered

$$P_1 = \frac{9}{s^2 + 1.2s + 9} \tag{51}$$

The LQR based pole placement tuning method, presented in section 2 as applied here with a closed loop parametric demand of $\zeta^{cl} = 0.75, \omega_n^{cl} = 7$ rad/sec; produces the optimal PID controller as in (52)

$$C_1 = 65.6944 + \frac{285.8333}{s} + 6.8667s \qquad (52)$$

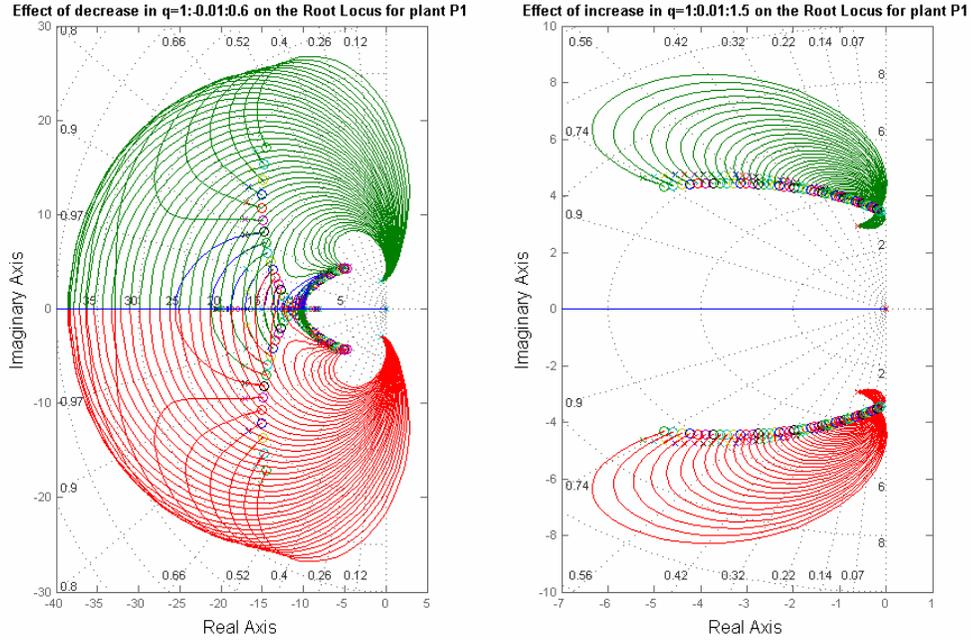

Fig. 4. Change in the position of the open loop controller zeros and dominant closed loop poles with gradual decrease and increase in $q$ for an under-damped process.

Now, with the variation in the controller differ-integral order $q$, of FOPID in (50), different set of controller parameters is obtained for the equivalent PID controller, representing the effect of FOPID controller (37) in the complex $s$-plane. The first observation is that the increase/decrease in the differ-integral order of FOPID ($q$), changes the position of the open loop controller zeros and closed loop system poles in the root locus-plot while location of the open loop system poles remain the same (Fig. 4). The second observation is that in each case, the dominant closed loop poles lie on the branch of the root locus starting from complex (conjugate) open loop system poles to open loop controller zeros depending on the system gain.

It is also evident from Fig. 4 that the increase in the order above $q = 1$ drastically reduces the closed loop damping indicating inferior time response. While, a gradual decrease in $q$ improves the closed loop damping but definitely up to a certain level. Once the open loop zeros of the derived sub-optimal PID touches the negative real axis, they move very fast towards lower damping region and therefore unstable region with further decrease in $q$, indicating highly oscillatory time response. Hence it is suggested to limit the closed loop pole locations before it touches the negative real axis with a desired closed loop damping slightly lesser than unity, by incorporating it as a constraint.

### 4.2. Critically-damped process

Next, a critically damped second order process is considered (53) having the open loop transfer function [13]

$$P_2 = \frac{25}{s^2 + 10s + 25} \tag{53}$$

The LQR based dominant pole placement tuning with desired closed loop performances as $\zeta^{cl} = 0.75, \omega_n^{cl} = 10\,\text{rad/sec}$ yields the following PID controller (54)

$$C_2 = 48 + \frac{300}{s} + 3.2s \tag{54}$$

The effect of variation in $q$ is similar as found in the previous subsection (for the lightly damped case) and is shown in Fig. 5.

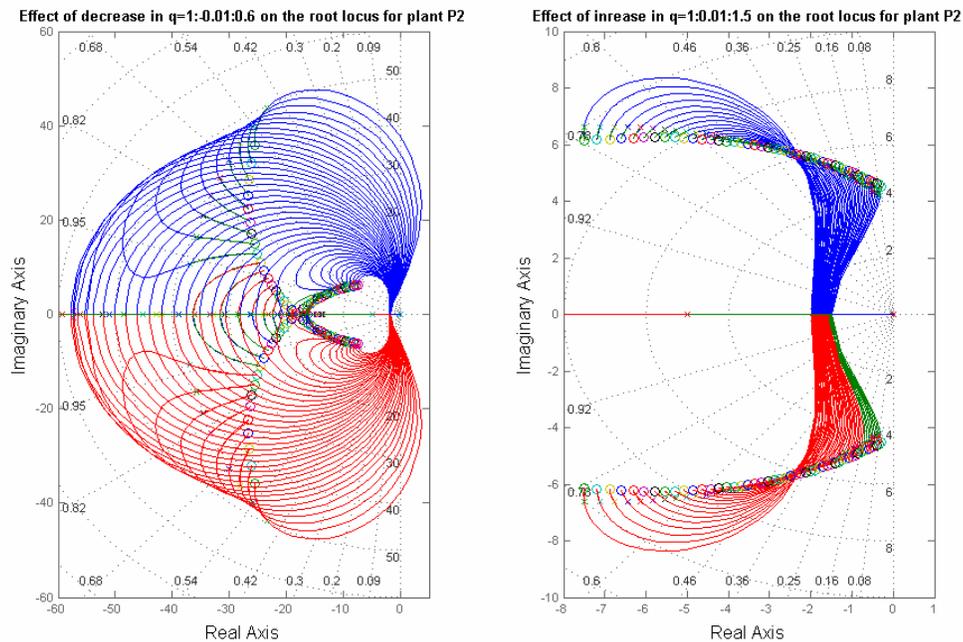

Fig. 5. Change in the position of the open loop controller zeros and dominant closed loop poles with gradual decrease and increase in $q$ for a critically-damped process.

### 4.3. Over-damped process:

A heavily damped second order process is (55) studied next which has the open loop transfer function [13]

$$P_3 = \frac{1}{s^2 + 10s + 1} \tag{55}$$

This over-damped process is tuned via LQR with a desired closed loop specification of $\zeta^{cl} = 0.75, \omega_n^{cl} = 5\,\text{rad/sec}$ which produces the optimal PID controller as in (56)

$$C_3 = 305.25 + \frac{937.5}{s} + 35s \tag{56}$$

The effect of variation in FOPID controller order ( $q$ ) on the derived PID controller parameters (51) is shown in Fig. 6 which is similar to the earlier observations.

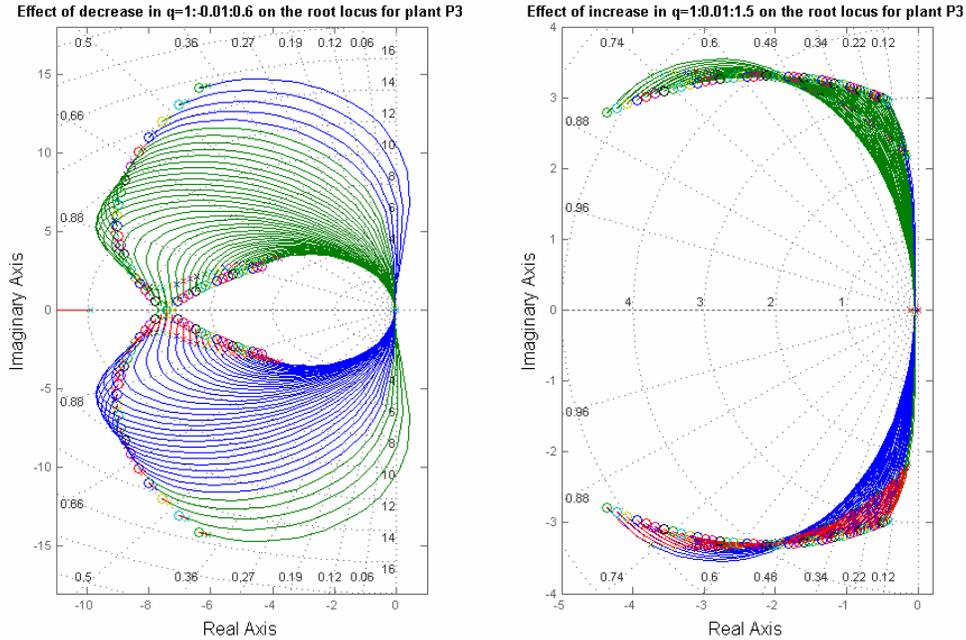

Fig. 6. Change in the position of the open loop controller zeros and dominant closed loop poles with gradual decrease and increase in $q$ for an under-damped process.

From the "M-curves" for the three different kind of second order processes as depicted in Fig. 4-6, it is evident that a judicious choice of $q$ is required to achieve desired closed loop performance. But the change of closed loop pole location with decrease in $q$ to place it at a desired damping could have been achieved with the dominant pole placement technique of PID controller itself and also in an optimal fashion with LQR. Then question arises what is the rationale to choose a sub-optimal method instead of LQR which is based on optimal pole placement at a comparatively lower damping as the first stage of tuning and then in the second stage to change the order $q$ until the closed loop pole reaches the desired damping? This issue is addressed in the following section with simulation examples which shows some advantage of this sub-optimal two stage PID tuning strategy; involving fractional order approach over the conventional LQR based pole-placement for PID controllers.

## 5. Two stage suboptimal tuning and its advantages over conventional LQR based tuning

In order to show the effectiveness of the proposed two stage sub-optimal PID controller tuning methodology, the basic design steps are given as follows.

***Step 1***:  For a given second order system of the form (2) with specified $\zeta^{ol}$ and $\omega_n^{ol}$ ,

**Step 1:** calculate the PID controller gains via pole placement for a low $\zeta^{cl}$ using equation (5).

**Step 2:** Location of the PID controller zeros can be intuitively modified using the approximated fractional order PID controller approach (50) by gradually decreasing the fractional order $q$, so as to meet desired $\zeta^{cl}$ (closer to a value which gives less overshoot).

**Step 3:** Using (5) obtain a single stage dominant pole placement at the same location as the desired $\zeta^{cl}$ and obtained $\omega_n^{cl}$ with the two stage tuning method.

**Step 4:** The cost of control for the dominant pole placement method (5) and two stage fractional order approach (50) can be obtained from the inverse optimal control formulation using (23) and (25) i.e. deriving the equivalent Riccati solutions from the each set of PID controller gains.

**Step 5:** Eigen-values of the differential Riccati solution (difference of the equivalent $P$ matrices) for PID controller gains in each case, indicate the comparative cost of control.

**Step 6:** Controller with low proportional-integral and derivative gains among the single stage and two stage tuning definitely indicates a lower control signal or controller effort and smaller actuator size with the time domain performance remaining same.

LQR based dominant pole placement tuning results for the three kinds of second order plants having the transfer functions (51), (53) and (55) are shown in section 4, with the obtained PID controllers (52), (54) and (56) respectively as the first stage of proposed PID controller tuning methodology. Now, one pole and two complex zeros of the PID controllers in the complex $s$-plane are replaced by the corresponding fractional order poles/zeros [9]-[10] while keeping the similar order of the poles and zeros with an extra tuning knob. This extra flexibility induced by the equivalent integer order PID controller zeros at a certain point in $s$-plane can be thought of as to get replaced by the corresponding fractional order zeros of a FOPID controller at the same position in the $s$-plane, suggesting preservation of the desired control action for specified closed loop damping and frequency. Using first two steps, presented in the tuning algorithm, the FOPID controller order $q$ is gradually decreased up to 0.9 which give the approximated sub-optimal PID controller that was initially tuned via LQR as in section 4. The damping and frequency of the dominant closed loop poles are calculated and presented in Table 1. The equivalent PID controller gains for the FOPID controllers are also reported in Table 1.

Table 1:
Closed loop performance and derived controller parameters for reduced value of $q$

| Plant | Order of the FO controller ($q$) | $\zeta^{cl}$ | $\omega_n^{cl}$ | $\widehat{K_p}$ | $\widehat{K_i}$ | $\widehat{K_d}$ |
|---|---|---|---|---|---|---|
| $P_1$ | 0.9 | 0.934 | 8.88 | 120.4848 | 535.8142 | 8.5059 |
| $P_2$ | 0.9 | 0.927 | 13.7 | 83.166 | 565.4015 | 3.6415 |
| $P_3$ | 0.9 | 0.914 | 6.31 | 619.9069 | 2005.4 | 51.9555 |

Table 2:
Optimal and sub-optimal controller gains and corresponding weighting matrices

| Plant | Controller | $K_p$ | $K_i$ | $K_d$ | $Q_1$ | $Q_2$ | $Q_3$ | $R$ |
|-------|-----------|-------|-------|-------|-------|-------|-------|-----|
| $P_1$ | LQR | 160.6263 | 726.6801 | 10.9252 | 528060 | 100500 | 86.5792 | 1.0 |
| | Suboptimal | 120.4848 | 535.8142 | 8.5059 | 287100 | 5499.5 | 47.8442 | 1.0 |
| $P_2$ | LQR | 135.5376 | 953.4577 | 5.696 | 909080 | 7017 | 26.1576 | 1.0 |
| | Suboptimal | 83.166 | 565.4015 | 3.6415 | 319680 | 2512.8 | 9.5204 | 1.0 |
| $P_3$ | LQR | 704.0603 | 2296.3 | 59.2081 | 5273100 | 179260 | 3281.6 | 1.0 |
| | Suboptimal | 619.9069 | 2005.4 | 51.9555 | 4021600 | 137030 | 2498.7 | 1.0 |

Table 3:
Riccati solution for the optimal and sub-optimal controllers

| Plant | Controller | $P_{11}$ | $P_{12}$ | $P_{13}$ | $P_{22}$ | $P_{23}$ | $P_{33}$ |
|-------|-----------|----------|----------|----------|----------|----------|----------|
| $P_1$ | LQR | 117450 | 8036 | 80.7422 | 1706.5 | 17.8474 | 1.2139 |
| | Suboptimal | 65093 | 4629 | 59.5349 | 989.8673 | 13.3872 | 0.9451 |
| $P_2$ | LQR | 130180 | 5812.2 | 38.1383 | 793.7882 | 5.4215 | 0.2278 |
| | Suboptimal | 47588 | 2285.1 | 22.6161 | 317.1408 | 3.3266 | 0.1457 |
| $P_3$ | LQR | 1619100 | 158920 | 2296.3 | 46490 | 704.0603 | 59.2081 |
| | Suboptimal | 1245200 | 124250 | 2005.4 | 36453 | 619.9069 | 51.9555 |

Here, it is attempted to place the closed loop poles at the same locations presented in Table 1, via LQR presented in section 2 and the corresponding controller gains and weighting matrices for the optimal controller is reported in Table 2. Therefore, the third step gives the single stage pole placement gains for the PID controller using (5) for handling the three different systems.

According to principle of an inverse regulator problem [5], [18]-[20], for a given set of stabilizing state feedback controller gains it is possible to find out a positive definite Riccati solution ($P$) and weighting matrices ($Q$ and $R$) for a given system. Since the controller gains obtained, by varying the controller order $q$ does not preserve the required optimality criterion (14), hence by formulating an inverse optimal control problem the obtained Riccati solution ($P$) will always be sub-optimal than that directly obtained from LQR. The Riccati solution ($P$) and weighting matrices ($Q$ and $R$) for the two stage sub-optimal inverse regulator and also the single stage optimal quadratic regulator based PID controllers have been reported in Table 3. In fact, the weighting matrices and the Riccati solutions for the sub-optimal and LQR based PID controller gains can be obtained using the inverse regulator formulation in (26), (23) and (25) respectively, as the fourth step of the design algorithm. As the fifth step, eigen-values of the differential Riccati solutions needs to be checked for comparison of the cost of control with the LQR based and proposed sub-optimal FO method.

From Table 3, it is evident that for all the above three processes, each element of the Riccati solution matrix ($P$) is higher for LQR based tuning which again indicates

$$x^T(0)P_{LQR}x(0) > x^T(0)P_{subopt}x(0) \qquad \Rightarrow J_{LQR} > J_{subopt} \qquad (57)$$

It is well established theory that for a set of fixed weighting matrices ($Q$ and $R$) there is no $P$ and corresponding state feedback gains which yields lesser cost of control than that with the LQR based technique. But in this specific case, the two stage suboptimal tuning method has much lower control cost than that of the LQR as the weighting matrices no longer remains constants and decreases with reduction in $q$. Indeed, the introduction of fractional poles and zeros in the complex $s$ plane makes the optimal regulator to have much lower control cost even lower than the cost of LQR based controllers. This typical behavior can be justified as the weighting matrix ($Q$) reduces in the equivalent integer order approximation, obtained by inverse regulator theory. Here, the word "sub-optimal" has been suggested for the fractional order approach of controller tuning as because calculation of the quadratic control cost using (34) is valid only for systems with integer order rational pole-zero models. Fractional order systems give fractional order state space formulation whose optimality criterion is different from that generally used in classical optimal control theory as reported in Agarwal [21]. Simulation examples shows that fractional order optimal regulators have lower control cost than that with the integer order ones which is analyzed in [22]-[25] in a detailed manner. Indeed, the concept of fractional pole-zero placement in the $s$-plane instead of conventional (integer) pole-zero placement perhaps may give some extra flexibility for optimal tuning of controllers which outperforms the classical integer order optimal state-feedback regulators [23].

The step input and load disturbance responses with the LQR and suboptimal techniques are shown in Fig. 7 along with the corresponding control signals. It is evident from Fig. 7 that the time-responses with the two tuning strategies are almost overlapping with high capability of suppressing load disturbances. But significant reduction in the control signal for two-stage suboptimal tuning is evident compared to LQR based tuning. The fact can be verified from PID controller gains itself in Table 2 with the two approaches, as also stated in the sixth step of the algorithm. This essentially implies that the approximation of the fractional zero positions in the complex $s$-plane with integer order zeros does not change the closed loop performance of a system in a significant manner, but only the cost of control and initial control effort have been found to get reduced when computed with an inverse optimal regulator formulation.

In Table 2, all the PID controller gains are higher for each cases of LQR based tuning than the suboptimal one. Now the control signal for the three cases can be evaluated with (18). It is also evident from the schematic in Fig. 2 that the initial value of the state variables i.e. $x(0)$ do not get affected by the controller gains and are indeed same for two different philosophies of PID controller design. If the two set of PID gains are represented by $K_{LQR}$ and $K_{subopt}$ respectively, then according to (18) we get (58)

$$\left[K_{LQR}x(0)\right] > \left[K_{subopt}x(0)\right] \quad \Rightarrow \quad \left[u_{LQR}(0)\right] > \left[u_{subopt}(0)\right] \tag{58}$$

The expression (58) shows that the proposed methodology enables us to design PID controller with less chance of actuator saturation, with lesser controller effort and associated cost involved due to large size of the actuator while also achieving the same demanded closed loop time response specifications.

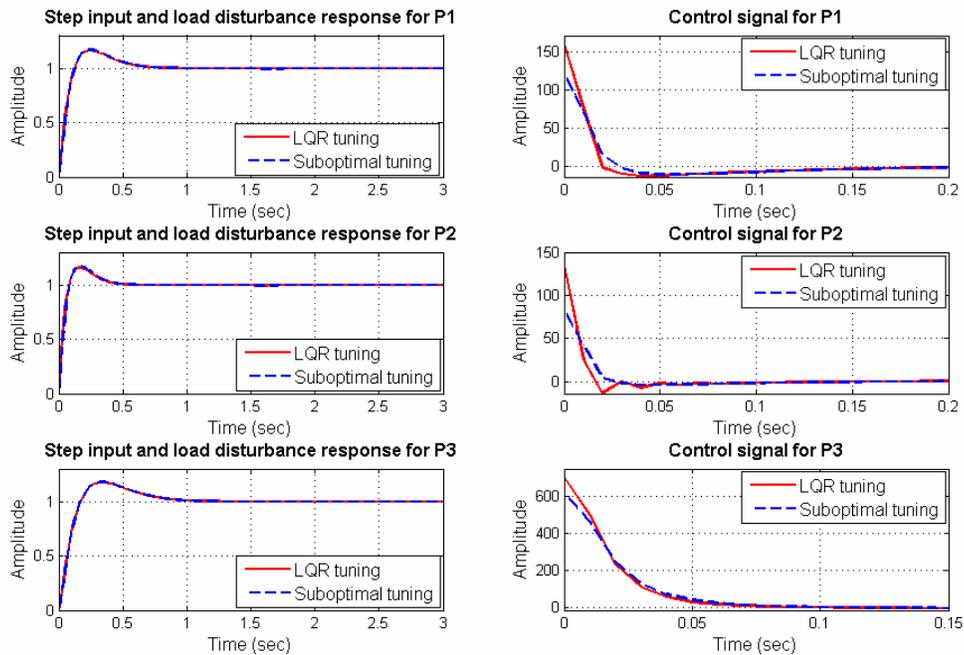

Fig. 7. Comparison of time responses and control signals for the test plants with LQR based and suboptimal PID tuning.

From the classical optimal control theory, it is already a well known fact that for a desired pole placement problem no methodology can give lesser control cost than that associated with an LQR formulation [15], [4]. But this conventional theory, breaks when fractional order poles and zeros comes into play as design variables. In the present study, it has been shown along with simulation examples that fractional poles and zeros ($q < 1$) of a FOPID controller, when approximated with integer order poles and zeros of an integer order PID in the same position of $s$-plane, produce lesser control cost and initial controller effort even lesser than that can be achieved from integer order quadratic optimal regulator, with no change in the closed loop response (Fig. 7).

## 6. Conclusion

A novel two stage sub-optimal PID controller tuning methodology is proposed in this paper, with the help of fractional order pole-zero placement approach. The methodology uses an LQR based optimal PID tuning approach with dominant pole placement in the first stage. Then considering the same closed loop pole location, achieved with a FOPID, the integro-differential orders of the FO controller is decreased in the second level of tuning to achieve the desired damping for the approximated integer order PID controller which exactly matches the same closed loop time domain and frequency domain characteristics where the design was initially attempted in a single shot. Simulation studies for three different classes of second order processes show that the control cost and most significantly the initial controller effort reduces to a large extent with the two stage sub-optimal technique over the LQR based one. Justification of finding better performance than the LQR based method is given using a fractional order concept of controller poles and zeros at the same location in the complex $s$-plane, over the LQR based guaranteed integer order pole placement. Performance comparison of the two tuning techniques is done while proposing a systematic computational methodology for the continuous Riccati solution and weighting matrices by the help of inverse optimal regulator theory. Step by step design technique is illustrated as a guideline for improved tuning of PID controller via fractional order pole-zero placements and their conformal mapping between $s \leftrightarrow w$ planes to achieve low cost of control and controller effort. The formulation proposed in this paper works well for second order (oscillatory or sluggish type) processes with negligibly small process delay. Further works could be extended towards finding similar optimal control formulation for delay dominant processes [26].

**Acknowledgement**

This work has been supported by the Board of Research in Nuclear Sciences (BRNS) of the Department of Atomic Energy, Govt. of India, sanction no. 2006/34/34-BRNS dated March 2007.

**Reference:**

[1] Karl J. Astrom and Tore Hagglund, "PID controller: Theory, design, and tuning", Instrument Society of America, 1995.

[2] Yu Zhang, Qing-Guo Wang, and K.J. Astrom, "Dominant pole placement for multi-loop control systems", *Automatica*, vol. 38, no. 7, pp. 1213-1220, July 2002.

[3] Qing-Guo Wang, Zhiping Zhang, Karl Johan Astrom and Lee See Chek, "Guaranteed dominant pole placement with PID controllers", *Journal of Process Control*, vol. 19, no. 2, pp. 349-352, Feb. 2009.

[4] Jian-Bo He, Qing-Guo Wang, and Tong-Heng Lee, "PI/PID controller tuning via LQR approach", *Chemical Engineering Science*, vol. 55, no. 13, pp. 2429-2439, July 2000.

[5] Brian D.O. Anderson and John B. Moore, "Optimal Control: linear quadratic methods", Prentice-Hall International, Inc., Englewood Cliffs, NJ, 1989.

[6] Igor Podlubny, "Fractional-order systems and $PI^{\lambda}D^{\mu}$-controllers", *IEEE Transactions on Automatic Control*, vol. 44, no. 1, pp. 208-214, Jan. 1999.

[7] Shantanu Das, "Functional fractional calculus", Springer, Berlin, 2011.

[8] Tom T. Hartley and Carl F. Lorenzo, "Dynamics and control of initialized fractional-order systems", *Nonlinear Dynamics*, vol. 29, no. 1-4, pp. 201-233, 2002.

[9] B.M. Vinagre, C.A. Monje, and A.J. Calderon, "Fractional order systems and fractional order control actions", *Lecture 3 IEEE CDC '02 TW#2: Fractional Calculus Applications in Automatic Control and Robotics*, 2002.

[10] Farshad Merrikh-Bayat, Mahdi Afshar, and Masoud Karimi-Ghartemani, "Extension of the root-locus method to a certain class of fractional-order systems", *ISA Transactions*, vol. 48, no. 1, pp. 48-53, Jan. 2009.

[11] F. Merrikh-Bayat and M. Karimi-Ghartemani, "Method for designing $PI^{\lambda}D^{\mu}$ stabilizers for minimum-phase fractional-order systems", *IET Control Theory and Applications*, vol. 4, no. 1, pp. 61-70, Jan. 2010.

[12] G.M. Malwatkar, S.H. Sonawane, and L.M. Waghmare, "Tuning PID controllers for higher-order oscillatory systems with improved performance", *ISA Transactions*, vol. 48, no. 3, pp. 347-353, July 2009.

[13] Rames C. Panda, Cheng-Ching Yu, and Hsiao-Ping Huang, "PID tuning rules for SOPDT systems: review and some new results", *ISA Transactions*, vol. 43, no. 2, pp. 283-295, April 2004.

[14] Mehrdad Saif, "Optimal linear regulator pole-placement by weight selection", *International Journal of Control*, vol. 50, no. 1, pp. 399-414, July 1989.

[15] Stanislaw H. Zak, "Systems and control", Oxford University Press, New York, 2003.

[16] Khalfa Bettou and Abdelfatah Charef, "Control quality enhancement using fractional $PI^{\lambda}D^{\mu}$ controller", *International Journal of Systems Science*, vol. 40, no. 8, pp. 875-888, August 2009.

[17] Riccardo Caponetto, Giovanni Dongola, Luigi Fortuna, and Antonio Gallo, "New results on the synthesis of FO-PID controllers", *Communications in Nonlinear Science and Numerical Simulation*, vol. 15, no. 4, pp. 997-1007, April 2010.

[18] T. Fujii, "A new approach to the LQ design from the viewpoint of the inverse regulator problem", *IEEE Transactions on Automatic Control*, vol. 32, no. 11, pp. 995-1004, Nov. 1987.

[19] Takao Fujii and Masaru Narazaki, "A complete solution to the inverse problem of optimal control", *21$^{st}$ IEEE Conference on Decision and Control, 1982*, vol. 21, pp. 289-294, Dec. 1982.

[20] P. Moylan and B. Anderson, "Nonlinear regulator theory and an inverse optimal control problem", *IEEE Transactions on Automatic Control*, vol. 18, no. 5, pp. 460-465, Oct. 1973.

[21] Om Prakash Agarwal, "A general formulation and solution scheme for fractional optimal control problems", *Nonlinear Dynamics*, vol. 38, no. 1-2, pp. 323-337, 2004.

[22] Christophe Tricaud and YangQuan Chen, "An approximate method for numerically solving fractional order optimal control problems of general form", *Computers and Mathematics with Applications*, vol. 59, no. 5, pp. 1644-1655, March 2010.

[23] Abdollah Shafieezadeh, Keri Ryan, and YangQuan Chen, "Fractional order filter enhanced LQR for seismic protection of civil structures", *ASME Journal of Computational and Nonlinear Dynamics*, vol. 3, no. 2, pp. 021404, April 2008.

[24] Yan Li and YangQuan Chen, "Fractional order linear quadratic regulator", *IEEE/ASME Conference on Mechatronic and Embedded Systems and Applications, MESA 2008*, pp. 363-368, Oct. 2008, Beijing.

[25] Om P. Agarwal, "A quadratic numerical scheme for fractional optimal control problems", *ASME Journal of Dynamic Systems, Measurement and Control*, vol. 130, no. 1, pp. 011010, Jan. 2008.

[26] Guo-Ping Cai, Jin-Zhi Huang, and Simon X. Yang, "An optimal control method for linear systems with time delay", *Computers & Structures*, vol. 81, no. 15, pp. 1539-1546, July 2003.